\documentclass[11pt,a4paper]{article}
\usepackage[english]{babel}
\usepackage{amsmath}
\usepackage{amssymb}

\newtheorem{theorem}{Theorem}

\newtheorem{corollary}[theorem]{Corollary}
\newtheorem{lemma}[theorem]{Lemma}

\newtheorem{remark}[theorem]{Remark}

\newcounter{tmpabcd}

\newcounter{tmpnum}

\newcounter{tmprome}

\newcommand{\mcc}{\mathbb{C}}

\newcommand{\mpp}{\mathbb{P}}

\newcommand{\mnn}{\mathbb{N}}

\newcommand{\mrr}{\mathbb{R}}

\newcommand{\mqq}{\mathbb{Q}}

\newcommand{\Dist}{\ensuremath{{\rm Dist}}}

\renewcommand{\b}[1]{{{#1}}}

\newcommand{\hidden}[1]{}

\newcommand{\ul}[1]{\underline{#1}}
\newcommand{\ull}[1]{\underline{\b{#1}}}

\newcommand{\ol}[1]{\overline{#1}}

\newcommand{\ord}{\ensuremath{{\rm ord}}}
\newcommand{\ordz}{\ensuremath{{\rm ord_{\b{z}=0}}}}

\newcommand{\trdeg}{\ensuremath{{\rm tr.deg.}}}

\begin{document}

\markboth{Evgeniy Zorin}
{Algebraic Independence and Mahler's method}

\title{Algebraic Independence and Mahler's method}

\author{EVGENIY ZORIN\footnote{Institut de math\'ematiques de Jussieu, Universit\'e Paris 7, Paris, France.} \footnote{E-mail: evgeniyzorin@yandex.ru}}
\date{}





\maketitle

\begin{abstract}
\selectlanguage{english}
We give some new results on algebraic independence within Mahler's method, including algebraic independence of values at transcendental points. We also give some new measures of algebraic independence for infinite series of numbers. In particular, our results furnishes, for $n\geq 1$ arbitrarily large, new examples of sets $\left(\theta_1,\dots,\theta_n\right)\in\mrr^n$ normal in the sense of definition formulated by Grigory Chudnovsky (1980).
\end{abstract}







\selectlanguage{english}

Let $p(\b{z})=p_1(\b{z})/p_2(\b{z})$ be a rational fraction with coefficients in $\ol{\mqq}$. We denote $d=\deg p$, $\delta=\ordz p$. Let $f_1(\b{z})$,...,$f_n(\b{z})\in\ol{\mqq}[[\b{z}]]$
be functions analytic in some neighborhood $U$ of $0$, with algebraic coefficients and satisfying the following system of functional equations
\begin{equation} \label{systeme_1}
    a(\b{z})\ull{f}(\b{z})=A(\b{z})\ull{f}(p(\b{z}))+B(\b{z}),
\end{equation}
 where $\ull{f}(\b{z})=(f_1(\b{z}),\dots,f_n(\b{z}))$, $a(\b{z})\in\ol{\mqq}[\b{z}]$ and $A$ (resp. $B$) is a matrix $n\times n$ (resp. $n\times 1$) with coefficients in $\ol{\mqq}[\b{z}]$.

Algebraic independence of values of such functions was studied by Becker, Mahler, Nishioka, T\"opfer and others~\cite{B1994,Ni1986,Ni1996,Pellarin2010,ThTopfer1995}. For this purpose one can also use a general method developed in~\cite{PP_KF} (see also~\cite{PP1997}). 
This method requires multiplicity estimate. Recently a new multiplicity lemma for solutions of~(\ref{systeme_1}) was established (see~\cite{EZ2010}, Theorem~3.11, also~\cite{EZ2011}). Using this latter result with the general method from~\cite{PP_KF} one can deduce the following theorems, which improve previously known and establish new results on algebraic independence and measures of algebraic independence.

\begin{theorem} \label{annexe_theo_ia1}
Let $f_1({z})$,\dots,$f_n({z})$ be analytic functions as described above, algebraically independent over $\mcc(z)$, and with additional condition $p({z})\in\ol{\mqq}[{z}]$. Let $y\in \ol{\mqq}^*$ be such that
\begin{equation*}
    p^{[h]}(y)\rightarrow 0
\end{equation*}
(as $h\rightarrow\infty$) and no iterate $p^{[h]}(y)\ne 0$ is a zero of $\det A$ .


Then for all $\varepsilon>0$ there is a constant $C>0$ such that for all variety $W\subset\mpp^n_{\mqq}$ of dimension $k < n+1-\frac{\log d}{\log\delta}$, one has
\begin{equation} \label{annexe_theo_ia1_result}
   \log\Dist(x,W)\geq-C\left(h(W)+d(W)^{\frac{n+1-k+\varepsilon}{n+1-k-\frac{\log
d}{\log\delta}}}\right)^{\frac{n+1}{n-k}-\frac{\log\delta}{\log d}\frac{k+1}{n-k}}
\kern-28pt
\times d(W)^{\frac{\log\delta}{\log d}\frac{k+1}{n-k}},
\end{equation}
where $x=\left(1:f_1(y):\dots:f_n(y)\right)\in\mpp^n_{\mcc}$.
\end{theorem}
\begin{remark}
The definition of $\Dist(x,W)$ (for a point $x\in\mpp^n$ and a subvariety $W$of the same space) can be found in \cite{NP}, Chapter~6, \S~5 or~\cite{EZ2010}, \S~1.3 (see Definition~1.17 {\it loc.cit.} and discussion after it)). Two simple cases are of crucial immportance to get the expression of this notion. First of all, if $W$ is a hypersurface defined by a homogeneous polynomial $P$, then $\Dist(x,W)$ is essentially $\|P(x)\|$ (more precisely, in this case $\log\Dist(x,W)=\log\|P(x)\|-\deg(P)\cdot\log\|x\|-\log\|P\|$). So essentially if $k=n-1$ is admitted in the statement of Theorem~\ref{annexe_theo_ia1}, one can substitute $\log\|P(x)\|$ in place of $\log\Dist(x,W)$ in the l.h.s. of~(\ref{annexe_theo_ia1_result}), obtaining the estimate which is usually considered itself as a measure of algebraic independence

On the other hand, for all points $x\in\mpp$ and all subvarieties $W\in\mpp^n$ one has $\Dist(x,W)=0$ iff $x\in W$. So if some value of $k$ is admitted in Theorem~\ref{annexe_theo_ia1} (i.e. if $k < n+1-\frac{\log d}{\log\delta}$), then at least $k+1$ of values $f_1(y),\dots,f_n(y)$ are algebraically independent over $\mqq$ (as the r.h.s. of~(\ref{annexe_theo_ia1_result}) is $>-\infty$). It immediately gives us two following corollaries.
\end{remark}
\begin{corollary} \label{annexe_cor_ia0}
In the situation of Theorem~\ref{annexe_theo_ia1} one has
\begin{equation*}
    \trdeg_{\mqq}\mqq\left(f_1(y),\dots,f_n(y)\right)\geq n+1-\left[\frac{\log d}{\log\delta}\right].
\end{equation*}
\end{corollary}
\begin{corollary} \label{annexe_cor_ia1}
In the situation of Theorem~\ref{annexe_theo_ia1} and assuming $\frac{\log d}{\log\delta}<2$, one has
\begin{equation} \label{annexe_cor_ia1_est} \index{Degre de transcendance, minoration@Degr\'e de transcendance, minoration}
    \trdeg_{\mqq}\mqq\left(f_1(y),\dots,f_n(y)\right)=n.
\end{equation}
\end{corollary}

Corollary~\ref{annexe_cor_ia0} improves on Theorem~3 of~\cite{ThTopfer1995} that gives a weaker above estimate:
\begin{equation*}
    \trdeg_{\mqq}\mqq\left(f_1(y),\dots,f_n(y)\right)\geq \lceil(n+1)\frac{\log\delta}{\log d}-1\rceil,
\end{equation*}
where $\lceil*\rceil$ denotes the smallest integer bigger than $*$.

Corollary~\ref{annexe_cor_ia1} improves on Corollary~2 of~\cite{ThTopfer1995} (the estimate~(\ref{annexe_cor_ia1_est}) was proven there only for $n=1$).

We can also give a measure of algebraic independence of values $f_1(y),\dots,f_n(y)$, in general case $y\in\mcc^*$. No result of this type was known before in this situation. Nevertheless, it should be mentioned that our estimates in this situation are weaker than in the case of algebraic $y$. 

\begin{theorem} \label{annexe_theo_ia2} Let $f_1({z})$,\dots,$f_n({z})$ be analytic functions as described in the beginning. We assume them to be algebraically independent over $\mcc(z)$ and we assume also $p({z})\in\ol{\mqq}[{z}]$ with $\delta=\ord_{z=0} p({z})\geq 2$  and $d=\deg p(z)$.
Let $y\in\mcc^*$ be such that
\begin{equation*}
    p^{[h]}(y)\rightarrow 0
\end{equation*}
(with $h\rightarrow\infty$) and no iterate $p^{[h]}(y)\ne 0$ is a zero of $\det A$.

Then for all $\varepsilon>0$ there is a constant $C$ such that for all variety $W\subset\mpp^{n+1}_{\mqq}$ of dimension $k < n+1-2\frac{\log d}{\log\delta}$, one has
\begin{equation*} \index{Mesure d'independance algebrique@Mesure d'ind\'ependance alg\'ebrique}
   \log\Dist(x,W)\geq-C\left(h(W)+d(W)^{\frac{n+1-k-\frac{\log
d}{\log\delta}+\varepsilon}{n+1-k-2\frac{\log
d}{\log\delta}}}\right)^{2\frac{n+1}{n-k}-\frac{\log\delta}{\log d}\frac{k+1}{n-k}}
\kern-42pt \times d(W)^{\frac{\log\delta}{\log d}\frac{k+1}{n-k}-\frac{n+1}{n-k}},
\end{equation*}
where $x=\left(1:y:f_1(y):\dots:f_n(y)\right)\in\mpp^{n+1}_{\mcc}$.
\end{theorem}
As before, one readily deduces two corollaries:
\begin{corollary} \label{annexe_cor_ia2}
In the situation of Theorem~\ref{annexe_theo_ia2} one has
\begin{equation*}
    \trdeg_{\mqq}\mqq\left(y,f_1(y),\dots,f_n(y)\right)\geq n+1-\left[2\frac{\log d}{\log\delta}\right].
\end{equation*}
\end{corollary}

\smallskip

\begin{corollary} \label{annexe_cor_ia3}
In the situation of Theorem~\ref{annexe_theo_ia2} and assuming $\frac{\log d}{\log\delta}<3/2$ one has
\begin{equation} \label{annexe_cor_ia3_est}
    \trdeg_{\mqq}\mqq\left(y,f_1(y),\dots,f_n(y)\right)\geq n-1.
\end{equation}
\end{corollary}

The next theorem improves Theorem~1 and Theorem~2 of~\cite{ThTopfer1995}, qualitatively et quantitatively.
\begin{theorem} \label{annexe_theo_2}
 Let $f_1({z})$,\dots,$f_n({z})$ be a set of functions analytic in some neighborhood of 0 and such as described in the beginning of this paper (in this statement we admit any $p(\b{z})\in\ol{\mqq}(\b{z})$ and not only $p(\b{z})\in\ol{\mqq}[\b{z}]$ as in preceding statements; we denote as always $d=\deg p$, $\delta=\ordz p\geq 2$). We suppose that these functions are algebraically independent over $\mcc(z)$. Assume $f_i(0)=0$, $i=1,\dots,n$, a number $y\in U\cap\ol{\mqq}$ satisfies $\lim_{m\rightarrow\infty}p^{[m]}(y)=0$ and for all $m\in\mnn$ a number $p^{[m]}(y)\ne 0$ is not a zero neither $\det A(\b{z})$  nor $a(\b{z})$. Then there is a constant $C>0$ such that for all variety $W\subset\mpp^n_{\mqq}$ of dimension $k < 2n+1-\frac{\log d}{\log\delta}(n+1)$, one has
\begin{equation*} \label{annexe_theo_2_result} \index{Mesure d'independance algebrique@Mesure d'ind\'ependance alg\'ebrique}
   \log\Dist(\ul{x},W)\geq-C\left(h(W)+d(W)^{\frac{1}{1-\frac{\log d}{\log\delta}\frac{n+1}{2n-k+1}}}\right)^{\frac{n+1}{n-k}-\frac{\log\delta}{\log d}\frac{k+1}{n-k}}d(W)^{\frac{k+1}{n-k}},
\end{equation*}
where $\ul{x}=\left(1:f_1(y):\dots:f_n(y)\right)\in\mpp^n_{\mcc}$.
In particular,
\begin{equation} \label{annexe_theo_2_trdeg} \index{Degre de transcendance, minoration@Degr\'e de transcendance, minoration}
    \trdeg_{\mqq}\mqq\left(f_1(y),\dots,f_n(y)\right)\geq 2n+1-\frac{\log d}{\log\delta}(n+1).
\end{equation}
\end{theorem}

Now we give a family of concrete examples with sets of functions arbitrarily long and satisfying all the hypothesis of our theorems. In the sequel we consider a particular case of~(\ref{systeme_1}) when this system is a diagonal one:
\begin{equation} \label{system_khi}
    \chi_i(z)=\chi_i\left(p(z)\right)+q_i(z), \quad i=1,\dots,n,
\end{equation}
where $p\in\ol{\mqq}(z)$ and $q_i\in\ol{\mqq}[z]$, $i=1,\dots,n$. Assuming $\deg q_i\geq 1$ and $q_i(0)=0$, $i=1,\dots,n$, $\ord_{z=0}p\geq 2$ we obtain solutions of~(\ref{system_khi}) analytic in some neighborhood of 0:
\begin{equation} \label{khi_explicit}
    \chi_i(z)=\chi_i\left(p(z)\right)+q_i(z), \quad i=1,\dots,n.
\end{equation}
We can deduce with Lemma~6 of~\cite{ThTopfer1995} as well as with Theorem~2 of~\cite{Kubota1977} the following lemma allowing to verify easily in many situations the algebraic independence of $\chi_1,\dots,\chi_n$ over $\mcc(z)$.
\begin{lemma} \label{annexe_chi_ind}
Let $n\in\mnn^*$, $q_i\in\mcc[\b{z}]$, $i=1,\dots,n$ and $p\in\mcc[\b{z}]$ satisfying $q_i(0)=0$, $i=1,\dots,n$, $p(0)=0$
and $p(\b{z})\ne\b{z}$. Let $\chi_1,\dots,\chi_n\in\mcc((\b{z}))$ be functions defined by~(\ref{khi_explicit}). Suppose that $1,q_1,\dots,q_n$ are $\mcc$-linearly independents and that at least one of the following conditions is satisfied:
\begin{enumerate}
  \item \label{annexe_chi_ind_condition_a} $\deg p \nmid \deg\left(\sum_{i=1}^n s_i q_i(\b{z})\right)$ for all $(s_1,\dots,s_n)\in\mcc^n\setminus\{\ul{0}\}$. \label{annexe_theo_khi_cond1}
  \item \label{annexe_chi_ind_condition_b} $\sum_{i=1}^n s_i \chi_i(\b{z}) \not\in \mcc[\b{z}]$ for all $(s_1,\dots,s_n)\in\mcc^n\setminus\{\ul{0}\}$. \label{annexe_theo_khi_cond2}
\end{enumerate}
Then the functions $\chi_1,\dots,\chi_n$ are algebraically independent over $\mcc(\b{z})$.
\end{lemma}
Using this lemma (especially point~(\ref{annexe_chi_ind_condition_a}) which is essentially due to Th.T\"opfer) we can assure in many situations algebraic independence of functions~(\ref{khi_explicit}) deducing immediately the conclusions of Theorems~\ref{annexe_theo_ia1}, \ref{annexe_theo_ia2}, \ref{annexe_theo_2} and their corollaries.

\begin{remark}
In~\cite{Ch1980} G.V.Chudnovsky introduced the notion of "normality" of $n$-uplets $(x_1,\dots,x_n)\in\mcc^n$. One says that  $(x_1,\dots,x_n)$ is \emph{normal} if it has a measure of algebraic independence of the form $\exp(-Ch(P)\psi(d(P)))$, i.e. if for all polynomial $P\in\ol{\mqq}[X_1,\dots,X_n]\setminus\{0\}$ one has the estimate
\begin{equation}\label{mesure_tau}
    |P(x_1,\dots,x_n)|\geq \exp(-Ch(P)\psi(d(P))),
\end{equation}
where $C>0$ is a real constant and $\psi:\mnn\rightarrow\mrr^{+}$ is an arbitrary function. Moreover, if one has the estimate~(\ref{mesure_tau}) with $\psi(d)=d^\tau$ for some constant $\tau$ one says that this $n$-uplet has a measure of algebraic independence of Dirichlet's type. In this situation one also defines Dirichlet's exponent to be the infimum of  $\tau$ admitted for $\psi(d)=d^{\tau}$ in~(\ref{mesure_tau}). In~\cite{Ch1980} G.V.Chudnovsky mentioned that for $n\geq 2$ the examples of normal $n$-uplets are quite rare, though almost all (in the sense of Lebesgue measure) $n$-uplets of complex numbers are normals. 

Th.T\"opfer gave a construction for a family of examples of normal $n$-uplets with Dirichlet's exponent $2n+2$ (see Theorem~1 and Corollary~4 of~\cite{ThTopfer1995}).

Our theorems 
assure the exponent of Dirichlet $n+2$ for a subfamily of these examples and allow also to produce new examples of normal $n$-uplets (due to the condition~(\ref{annexe_chi_ind_condition_b}) of Lemma~\ref{annexe_chi_ind}).
\end{remark}

\begin{center}%
          {\bfseries Acknowledgement\vspace{-.5em}}%
\end{center}%
The author would like to express his profound gratitude to Patrice \textsc{Philippon}. His interventions at many stages of  this research was of decisive importance. 

{\small



\begin{thebibliography}{0}

\bibitem{B1994} P.-G.Becker, "Transcendence measures for the values of generalized Mahler functions in arbitrary
characteristic", Publ. Math. Debrecen 45 (1994), 269-282.

\bibitem{Ch1980} G.V.Chudnovsky, "Measures of irrationality, transcendence and algebraic independence. Recent progress", Journ\'ees Arithm\'etiques~1980 (J.Armitage, ed.), Cambridge Univ. Press, 1982, 11-82.

\bibitem{Kubota1977} K.K.Kubota, "On the algebraic independence of holomorphic solutions of certain functional equations and their values", Math.Ann. 227 (1977), 9-50.

\bibitem{NP} Yu.Nesterenko, Patrice Philippon (eds.),
 "Introduction to Algebraic Independence Theory", Vol.~1752, 2001, Springer.

\bibitem{Ni1986} K.Nishioka, "Algebraic independence of certain power series of algebraic numbers", J. Number Theory 23 (1986), 353-364.

\bibitem{Ni1996} K.Nishioka, "Mahler Functions and Transcendence", Lecture Notes in Math. 1631, Springer,
1996.

\bibitem{Pellarin2010} F.Pellarin, "An introduction to Mahler's method for
transcendence and algebraic independence", preprint, 2010. Disponible at http://hal.archives-ouvertes.fr/hal-00481912/fr/

\bibitem{PP1997} P.Philippon, "Une approche m\'ethodique pour la transcendance et l'ind\'ependance alg\'ebrique de valeurs de fonctions analytiques".
J. Number Theory { 64} (1997) 291-338.

\bibitem{PP_KF} P.Philippon, "Ind\'ependance alg\'ebrique et $K$-fonctions". J. reine angew. Math. 497 (1998), 1-15.

\bibitem{ThTopfer1995} Th.T\"opfer, "Algebraic independence of the values of generalized Mahler functions", Acta
Arithmetica, LXX.2 (1995).

\bibitem{EZ2010} E.Zorin, "Lemmes de z\'eros et relations fonctionnelles", th\`ese de doctorat de l'Universit\'e Paris 6, 2010. Accessible at http://tel.archives-ouvertes.fr/tel-00558073/fr/

\bibitem{EZ2011} E.Zorin, "Zero Order Estimates for Analytic Functions", pr\'epublication, arXiv: 1103.1174.

\end{thebibliography}

\def\cprime{$'$} \def\cprime{$'$} \def\cprime{$'$} \def\cprime{$'$}
  \def\cprime{$'$} \def\cprime{$'$} \def\cprime{$'$} \def\cprime{$'$}
  \def\cprime{$'$} \def\cprime{$'$} \def\cprime{$'$} \def\cprime{$'$}
  \def\cprime{$'$} \def\cprime{$'$} \def\cprime{$'$} \def\cprime{$'$}
  \def\polhk#1{\setbox0=\hbox{#1}{\ooalign{\hidewidth
  \lower1.5ex\hbox{`}\hidewidth\crcr\unhbox0}}} \def\cprime{$'$}
  \def\cprime{$'$}

\bigskip

\noindent{\footnotesize EZ\,: }\begin{minipage}[t]{0.9\textwidth}
\footnotesize{\sc Institut de math\'ematiques de Jussieu, Universit\'e Paris 7, Paris, France}\\
{\it E-mail address}\,:~~ \verb|zorin@math.jussieu.fr|\quad\emph{or}\quad\verb|EvgeniyZorin@yandex.ru|
\end{minipage}

}

\end{document}